\documentclass[a4paper]{amsart}

\usepackage{amsmath,amssymb,amsthm}

\usepackage{tikz}
\usetikzlibrary{shapes}
\usetikzlibrary{intersections}

\usepackage{pgfplots}
\usepackage{caption}
\usepackage{subcaption}

\usepackage{hyperref,caption}
\usepackage{graphicx}
\usepackage{graphics}

\newtheorem*{thm}{Theorem}
\newtheorem*{lemma}{Lemma}

\newtheorem{definition}{Definition}

\newcommand{\dist}{\operatorname{dist}}

\begin{document}

\title[]{Detecting localized eigenstates \\of linear operators}

\thanks{The research of J.L. was supported in part by the National Science Foundation under award DMS-1454939. He would also like to thank Yingzhou Li for helpful discussions.}

\author[]{Jianfeng Lu}
\address[Jianfeng Lu]{Department of Mathematics, Department of Physics, and Department of Chemistry,
Duke University, Box 90320, Durham NC 27708, USA}
\email{jianfeng@math.duke.edu}

\author[]{Stefan Steinerberger}
\address[Stefan Steinerberger]{Department of Mathematics, Yale University, New Haven CT 06511, USA} \email{stefan.steinerberger@yale.edu}

\begin{abstract} We describe a way of detecting the location of localized eigenvectors of a linear
system $Ax = \lambda x$ for eigenvalues $\lambda$ with $|\lambda|$ comparatively large.  We define the family of functions
$f_{\alpha}: \left\{1.2. \dots, n\right\} \rightarrow \mathbb{R}_{}$
$$ f_{\alpha}(k) = \log \left( \| A^{\alpha} e_k \|_{\ell^2} \right),$$
where $\alpha \geq 0$ is a parameter 
and $e_k = (0,0,\dots, 0,1,0, \dots, 0)$ is the $k-$th standard basis
vector.  We prove that eigenvectors associated to eigenvalues with
large absolute value localize around local maxima of $f_{\alpha}$: the
metastable states in the power iteration method (slowing down its
convergence) can be used to predict localization. We present a
fast randomized algorithm and discuss different examples: a random
band matrix, discretizations of the local operator $-\Delta + V$ and
the nonlocal operator $(-\Delta)^{3/4} + V$.
\end{abstract}
\maketitle

\section{Introduction and Main Idea}
\subsection{Introduction} We are interested in spatially localized
eigenvectors of matrices $A \in \mathbb{R}^{n \times n}$. These
objects are of paramount importance in many fields of mathematics: the
ground state and low-frequency behavior of quantum systems
\cite{anderson, elilu, lots, fil, Marzari:12}, the behavior of
metastable random dynamical systems \cite{bov, bouv1, bouv2}, the
detection of central points in graphs clusters \cite{cheng1}, the
principal component analysis for sample covariance matrix \cite{zht},
and many more.\\
The purpose of this paper is to introduce a simple idea,
which provably detects localized eigenstates associated to eigenvalues
with large absolute value at low computational cost. We introduce the
entire relevant theory for matrices $A \in \mathbb{R}^{n \times n}$,
however, a crucial ingredient is the following: when numerically
computing solutions for many infinite-dimensional linear operators of
interest (linear/nonlocal/fractional partial differential equations,
integral equations, ...), these are usually discretized and the
discretization respects the spatial ordering of the underlying
domain. In particular, if the original continuous object has localized
eigenstates and the discretization is sufficiently accurate, then the
discretized linear operator will have localized eigenstates on the
associated graph.  We will completely ignore the question of how
operators are discretized and restrict ourselves to the question of
how to find localized eigenvectors.

\subsection{Main idea}\label{sec:mainidea}

We are given a matrix $A \in \mathbb{R}^{n \times n}$ (not necessarily
symmetric) and are interested in finding, if they exist, the location
of localized eigenvectors concentrating their mass on relatively few
coordinates of
$$ Ax = \lambda x$$
for $\lambda$ in the spectral edge (meaning that $|\lambda|$ is
comparatively large to the rest of the spectrum, the low-lying eigenvalues, $|\lambda|$ close to 0, can be obtained via the very same method after a transformation of $A$, see below).   Since strongly localized eigenstates
are essentially created by localized structure, they should also be
detectable by completely local operations.

\begin{definition}
  We define
  $f_{\alpha}: \left\{1, 2, \ldots, n\right\} \rightarrow \mathbb{R}_{}$
  given by
  \begin{equation*}
    f_{\alpha}(k) = \log \left( \| A^{\alpha} e_k \|_{\ell^2} \right),
  \end{equation*}
  where $\alpha \geq 0$ is a parameter and
  $e_k = (\delta_{ik})_{i=1}^{n}$ is the $k-$th standard basis vector.
\end{definition}

The main idea is rather simple: if highly localized eigenvectors exist, then they have a nontrivial inner product with one of the standard basis vectors (whose size can be bounded from below depending only on the scale of localization and not on $n$). An iterated application of the matrix will then lead to larger growth than it would in other regions. The idea is vaguely related to the stochastic interpretation \cite{stein} of the Filoche-Mayboroda landscape function \cite{fil}. The logarithm counteracts the exponential growth purely for the purpose of visualization.

\medskip 

\begin{figure}[ht]
\centering
\includegraphics[width = 0.65\textwidth]{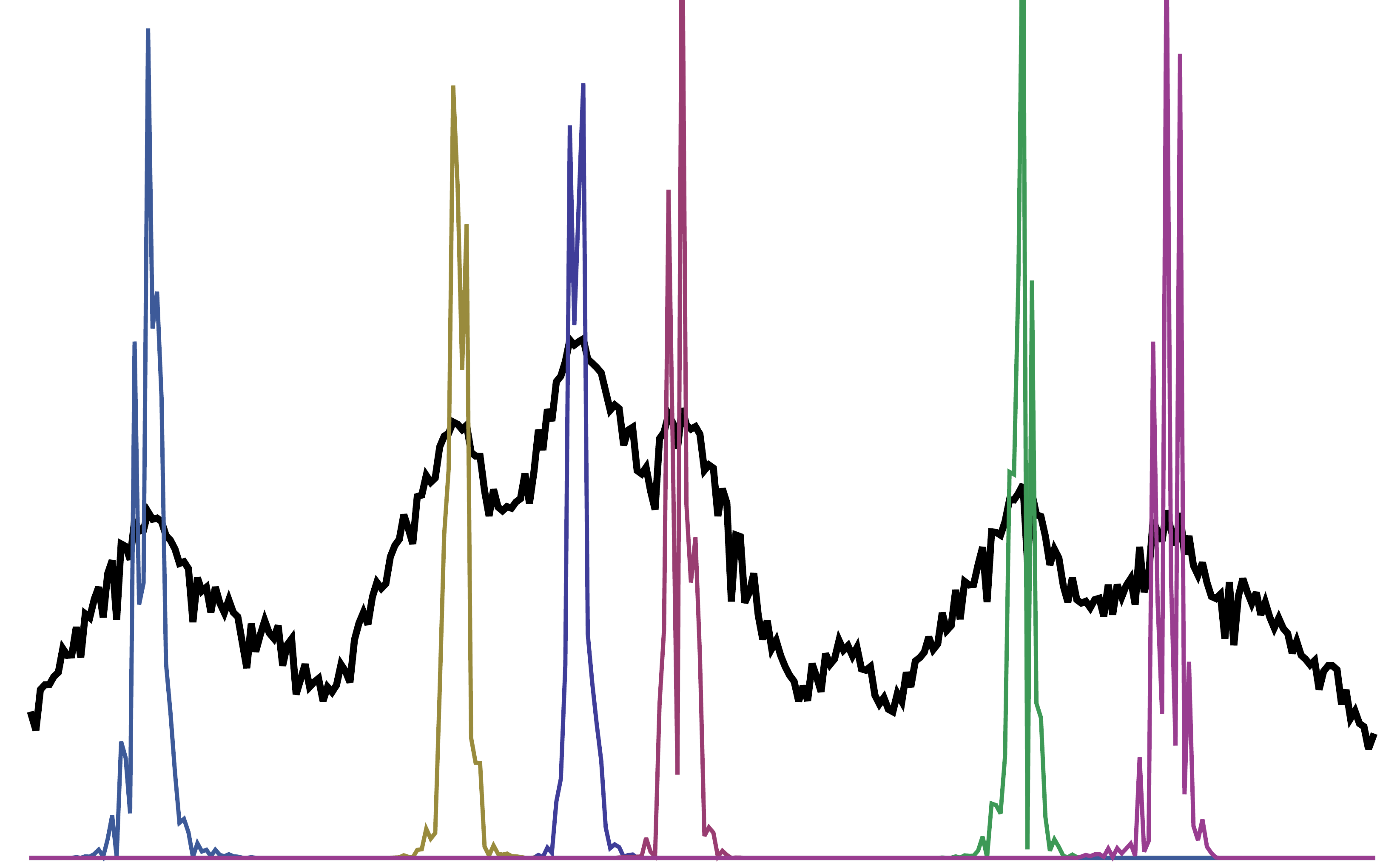}
\caption{$f_{50}$ (black) and the plots of the absolute value of the
  first $6$ eigenvectors associated to the $6$ largest eigenvalues.}
\label{fig:eins}
\end{figure}

\textit{Example.} We start by considering a numerical example (see
Fig. \ref{fig:eins}): here $A = B + B^{\top}$ where
$B \in \mathbb{R}^{n \times n}$, $n=300$, is given by a random band
matrix with bandwidth $2$ around the diagonal (i.e., $A$ is a $5-$diagonal
matrix) and every non-zero entry chosen independently randomly in the
interval $[-1,1]$. A typical outcome can be seen in
Figure~\ref{fig:eins} for $\alpha = 50$: the function $f_{\alpha}$ has
a series of local maxima and the first few eigenvectors localized
around these maxima; higher peaks in the landscape corresponds to
eigenvalues with larger absolute value.

The value of $\alpha$ depends on the precise circumstances; larger
values can lead to higher accuracy but also increase the computational
cost. It is worth pointing out that it is not interesting to have
$\alpha$ very large: whenever the largest eigenvalue $\lambda_1$ is
simple and the associated eigenvector $\phi_1$ does not vanish, then
this approach becomes less effective since
\begin{equation*}
  \lim_{\alpha \rightarrow \infty}{ \frac{ f_{\alpha}(k)}{\alpha}} = \bigl\lvert \bigl\langle \lambda_1 \phi_1, e_k\bigr\rangle\bigr\rvert.
\end{equation*}
It is not difficult to see that the convergence speed is going to depend on the spectral gap between the largest eigenvalue
and the rest (in terms of absolute value of the eigenvalues)-- it is commonly desirable to have a spectral gap; here, we are bound to encounter a delicate interplay between
spectral gaps and the scale of localization.

\section{Statament of Main result} 
We give one of the many possible formulation of a rigorous guarantee
of the approach. Indeed, the underlying principle, as outlined in the
previous section, is so simple that there are many ways of turning it
into a precise statement; we give a fairly canonical one but it is by
no means unique and different circumstances may call for different
versions.\\

\begin{figure}[h!]
\centering 
\begin{tikzpicture}[scale=0.8]
\draw [ultra thick] (-2,0) -- (10,0);
\filldraw (-2,0) circle (0.08cm);
\node at (-2, -0.4) {$1$};
\filldraw (10,0) circle (0.08cm);
\node at (10, -0.4) {$n$};
\draw [thick] (0,-0.2) -- (0,0.2);
 \draw[scale=1,thick, samples=50,domain=-2:10,smooth,variable=\x]  plot ({\x},{2*exp(-4*(\x-1)*(\x-1))});
 \node at (1, -0.4) {$J_1$};
 \draw [thick] (2,-0.2) -- (2,0.2);
 \draw [thick] (4,-0.2) -- (4,0.2);
 \draw[scale=1,thick, samples=50,domain=0:10,smooth,variable=\x]  plot ({\x},{2*exp(-4*(\x-3)*(\x-3))});
 \node at (3, -0.4) {$J_3$};
 \draw [thick] (5,-0.2) -- (5,0.2);
 \draw [thick] (6,-0.2) -- (6,0.2);
 \draw [thick] (7,-0.2) -- (7,0.2);
 \draw [thick] (8,-0.2) -- (8,0.2);
 \draw[scale=1,thick, samples=50,domain=0:10,smooth,variable=\x]  plot ({\x},{2*exp(-4*(\x-6)*(\x-6))});
 \node at (6, -0.4) {$J_4$};
 \draw[scale=1,thick, samples=50,domain=0:10,smooth,variable=\x]  plot ({\x},{2*exp(-4*(\x-7)*(\x-7))});
 \node at (7, -0.4) {$J_2$};
\end{tikzpicture}
\caption{The modulus square of four eigenvectors localized in four intervals $J_1, J_2, J_3, J_4$ bounded uniformly in length
$J_i \leq \mathcal{J}$. The eigenvectors have exponential decay away from the intervals. \label{fig:localeigenfun}}
\end{figure}
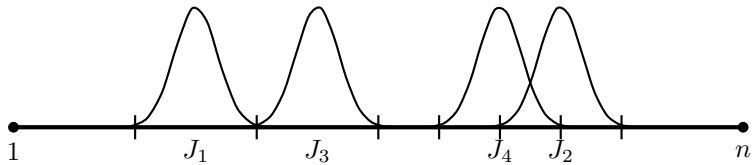

We start by clarifying our setup and introducing
some parameters below. We first phrase everything in a way that is
most natural in the setting of band matrices or matrices with rapid
decay off the diagonal (which covers $-\Delta + V$ on subsets of
$\mathbb{R}$, discretized by finite difference or finite volume
methods) -- the general case follows in a rather straightforward
manner by replacing the notion of `interval' by `subset', we briefly
discuss this below.  The only restrictive assumption is
the orthogonality of eigenvectors (1), which is usually given in the
setting that we are interested in (localization of self-adjoint
operators).  (2) and (3) introduces various parameters that are always
defined, however, in the non-localized regime they may result in a
vacuous conclusion (see Figure~\ref{fig:localeigenfun} for an illustration).
\begin{enumerate}
\item The eigenvectors of $A \in \mathbb{R}^{n \times n}$ form an orthonormal basis of $\mathbb{R}^n$ and we order the eigenvalues via
$$ |\lambda_1| \geq |\lambda_2| \geq |\lambda_3| \geq \dots \geq |\lambda_n|.$$
\item Every one of the first $k$ eigenvectors $\phi_i$ has half of its $\ell^2-$mass supported on an interval $J_{i} \subset \left\{1, 2, \dots, n \right\}$, i.e.
$$ \sum_{j \in J_i}{\phi_i(j)^2} \geq \frac{1}{2}$$
and we define $\mathcal{J}$ as the longest such interval
$$ \mathcal{J} = \max_{1 \leq i \leq k}{ |J_i|}.$$
\item We assume that, for all $1 \leq i \leq k$, the eigenvector $\phi_i$ has exponential decay away from the interval $J_i$, i.e. for all $1 \leq i \leq k$
$$ |\left\langle \phi_i, e_m \right\rangle| \leq \exp\left( - \beta \cdot \dist(m, J_i)\right)$$
for some universal constant $\beta > 0$. 
\end{enumerate}

Our main result states that, depending on the size of the spectral gap
and the quality of localization, there exist
$\alpha, y \in \mathbb{R}$ such that the superlevel set
$ \left\{ x: f_{\alpha}(x) \geq y\right\}$ intersects \textit{all}
localized intervals $J_i$ and $f_{\alpha}$ can only be that large in a
small neighborhood of the intervals $J_i$. In particular, this allows
for a detection of localization of the first $k$ eigenfunctions by
looking at $f_{\alpha}$ alone.

\begin{thm} If $\alpha \in \mathbb{N}$ is chosen such that 
\begin{equation*}
\left( \frac{\lvert \lambda_k \rvert}{\lvert \lambda_{k+1} \rvert} \right)^{\alpha} \geq 16 \sqrt{\mathcal{J}n},
\end{equation*}
then there exists a critical value $y \in \mathbb{R}$ such that $f_{\alpha}$ is large on all $J_i$
$$ \min_{1 \leq i \leq k}{ \max_{x \in J_i}{ f_{\alpha}(x)}} \geq y$$
and only large in their neighborhood: if $f_{\alpha}(m) \geq y - 1$, then
$$ \dist\left(m, \bigcup_{i=1}^{k}{J_i}\right)   \leq \frac{1}{ \beta} \left(\log{\left(16\sqrt{k}\sqrt{\mathcal{J}}\right)} + \alpha \log{\left( \left| \frac{\lambda_1}{\lambda_k} \right| \right)}\right).$$
\end{thm}

The condition on $\alpha$ depends on the spectral gap and localization
properties. If the spectral gap is large, then $\alpha \sim 1$ is
sufficient.  If the matrix satisfies that each row has a constant
number of non-zero entries (for example, if it is a local
discretization of a differential operator in $d$ dimensions), this
implies the computation of $f_{\alpha}$ may only require
$\mathcal{O}(n)$ operations.  It is clear from the proof that there
are many other possible conditions under which similar result could be
obtained. Natural variations include the following:
\begin{itemize}
\item The statement guarantees a gap of size $\geq 1$ between the
  values of $f_{\alpha}$ attained on the intervals and far away from
  their supports; in most instances a much smaller gap would suffice
  (especially if one identifies regions of localization via a notion
  of local maximum of $f_{\alpha}$, i.e., attaining a larger value
  than in a neighborhood of a certain size). We observe that our
  approach easily implies the inequalities
  $$ \alpha \log{|\lambda_1|} - \frac{1}{2} \log{n} \leq \max_{1 \leq m \leq n}{f_{\alpha}(m)} \leq \alpha \log{|\lambda_1|}$$
  and for any $i$ 
  $$  \max_{m \in J_i}{f_{\alpha}(m)} \geq \alpha \log{|\lambda_i|} - \frac{1}{2}\log{n}.$$
  In combination, they suggest that the gap size $\sim 1$ could in
  various situations be replaced by something much smaller which would
  improve the error bounds.
\item Orthogonality is not crucial: if the eigenvectors of $A$ form a
  basis and the angles between different eigenvectors are not too
  small, our argument easily implies similar results.
\item If we know in advance that a generic localized eigenfunction is going to be roughly localized on an interval
 $J$  then it is clear that it suffices to compute $f_{\alpha}$ on a
  $\sim |\mathcal{J}|-$net which further speeds up computation time.
\item $f_{\alpha}$ is defined via purely local operations and thus, by
  definition, its value at a certain location is stable under
  perturbing the matrix entries far away. On the other hand, it is
  entirely possible that a large perturbation will destroy the
  spectral structure of the matrix. With additional assumptions on the
  perturbation to guarantee spectral stability, it is then possible to
  use $f_{\alpha}$ of the unperturbed matrix to predict localization
  after perturbation. 
\end{itemize}

We also remark that the assumptions of the theorem consists of both
the existence of a gap of the spectrum and localization of the first
few eigenvectors. In some situations, the localization of the
eigenvectors follows from the gap assumption alone, for instance when
the matrix $A$ comes from a local discretization of a differential operator,
as established in \cite{Benzi, LinLu}.\\

As $\alpha$ increases, the neighborhoods in which the result guarantees localization are growing linearly in size (though only very slowly if $|\lambda_1/\lambda_k| \sim 1$) and $f_{\alpha}(m) \geq y-1$ becomes less informative. This is not an artifact but necessary:  whenever the largest eigenvalue $\lambda_1$ is simple and the associated
eigenvector $\phi_1$ never vanishes, then
$$ \lim_{\alpha \rightarrow \infty}{ \frac{ f_{\alpha}(k)}{\alpha}} = |\left\langle \lambda_1 \phi_1, \left( \delta_{ik}\right)_{i=1}^{n}\right\rangle|.$$
This is similar in spirit to the classical power method for computing the first eigenvector. 
It is noteworthy that we exploit exactly the fact that $|\lambda_i/\lambda_j| \sim 1 \pm \varepsilon$ in the edge of the spectrum which makes the power method a slow
method in practice. Or, put differently, our method exploits that highly localized eigenvectors associated to eigenvalues in the spectral edge \textit{correspond to metastable states for the power iteration!}
When $\lvert \lambda_1 \rvert \gg \lvert \lambda_k \rvert$, the first
eigenvector (and potentially other high lying ones) will then
interfere with the performance of the procedure. In such situations, we may
revise the procedure by first identifying those dominant eigenvectors and
then applying the procedure while iteratively projecting onto the orthogonal
complement of the subspace spanned by the dominant eigenvectors.
 The details are standard and left to the interested reader.\\

 In many applications, we have the eigenvalue problem 
$$ Ax = \lambda x$$
with $A$ being positive definite and the dominant characteristics of the physical system being determined by the low-lying eigenvalues $\lambda$. A straightforward application of
our localization technique
is only going to yield the largest eigenvalues. The obvious modification is to consider the matrix 
$$ \mbox{Id}_{n \times n} - \frac{A}{\|A\|}  \quad \mbox{instead.}$$
This operation preserves sparsity and flips the spectrum and the low-lying eigenvalues are now in the spectral edge. This is used in our numerical examples in Section~\ref{sec:num}.
Alternatively, if we are given a self-adjoint, positive definite and linear map on a Hilbert space $A:H \rightarrow H$ such that $A^{-1}$ is compact, then a natural way of recovering the bottom of the spectrum is
via considertion of the semigroup
$$ u_t + Au = 0.$$
An application of the spectral theorem allows us to write the semigroup as
$$ e^{-tA}u= \sum_{k \geq 0}{ e^{-\lambda_k t} \left\langle \phi_k, u \right\rangle \phi_k},$$
which has slow decay for the small eigenvalues and large decay for
larger eigenvalues. Note that if we Taylor expand $e^{-tA}$ and keeps
only the leading order term, we get $e^{-tA} \approx \mbox{Id} - t A$,
which connects with the previous trick. We refer to the Appendix A for an application.\\

 The result, as given above, is easiest to understand in the setting of
banded matrices. Banded matrices correspond naturally to localized interactions on the lattice
$\mathbb{Z}$. Neighborhoods of points correspond to intervals and this is how the Theorem
was phrased. At a greater level of generality, there need not be such underlying structure and
we will replace intervals by general subsets $J_i \subset \left\{1,2,\dots, n\right\}$. The notion
of distance between an element $k \in \left\{1,2, \dots, n\right\}$ and a subset $J_i$ is implicitly defined
via assuming the inequality
$$ |\left\langle \phi_i, e_k \right\rangle| \leq \exp\left( - \beta \cdot \dist(k, J_i)\right)$$
to hold. We emphasize that in all the interesting applications, where $A$ is the discretization of
a differential operator (or somewhat localized integral operator), these notions can be made
rather precise and we recover the classical notion of distance in Euclidean space.

\section{Fast Randomized Algorithms}
Computationally, $f_{\alpha}$ can be obtained by
calculating the $\ell^2$-norm of the rows of the matrix
$A^{\alpha}$. Thus the algorithm is particularly efficient if $A$ has
structures enabling fast multiplications, such as being sparse, low
rank, etc.  To further accelerate the computation, we exploit the
ideas from randomized numerical linear algebra (see \cite{Halko:11}
for a review) to use the following randomized version of the landscape
function. For simplicity, we assume that the matrix $A$ is symmetric in this section. 

\begin{definition} We define $f_{R, \alpha}: \left\{1, 2, \ldots, n\right\} \rightarrow \mathbb{R}$
  given by
  \begin{equation*}
    f_{R, \alpha}(k) = \log \bigl( \lVert e_k^{\top} A^{\alpha} R \rVert_{\ell^2} \bigr), 
  \end{equation*}
  where $\alpha \geq 0$ is a parameter,
  $e_k = (\delta_{ik})_{i=1}^{n}$ is the $k-$th standard basis vector,
  and $R$ is a $n \times m$ random matrix with i.i.d. $N(0, 1)$ entries.
\end{definition}

In terms of computation, if $m \ll n$, the randomized version
$f_{R, \alpha}$ only requires applying the matrix $A$ on a tall skinny
matrix for $\alpha$ times, and thus, the randomized version is
particularly advantageous for dense $A$, as it brings down the cost
from $\mathcal{O}(n^3)$ to $\mathcal{O}(n^2 m)$. As it turns out, the efficiency of this
method is intimately coupled to very well studied concepts centered
around the stability of random projection onto subspaces. To see this,
we denote the columns of $R$ by $r_1, \dots, r_m$ and observe that
$$  \lVert e_k^{\top} A^{\alpha} R \rVert_{\ell^2}^2 = \sum_{i=1}^{n}{ \bigl\lvert \lambda_i^{\alpha} \phi_{i, k}\bigr\rvert^2 \sum_{j=1}^{m}{ |\left\langle \phi_i, r_j \right\rangle|^2 }},$$
while for $f_{\alpha}$, we have 
\begin{equation*}
  \lVert e_k^{\top} A^{\alpha} \rVert^2 = \sum_{i=1}^n \bigl\lvert \lambda_i^{\alpha} \phi_{i, k}\bigr\rvert^2.
\end{equation*}
We start by quickly discussing a very strong sufficient condition that
allows to transfer results almost verbatim from our Theorem to the
random case and is equivalent to classical questions in dimensionality
reduction; stronger results are discussed below. We observe that
certainly all the results transfer in a pointwise manner if we knew
that
$$ \sum_{j=1}^{m}{ |\left\langle \phi_i, r_j \right\rangle|^2 } \qquad \mbox{is almost constant in the parameter}~i$$
for a typical realization of $m$ such random vectors.  However, since
the $\phi_i$ span the space, we are really asking that the map
$h:\mathbb{S}^{n-1} \rightarrow \mathbb{R}$ given by
$$ h(v) = \sum_{j=1}^{m}{ |\left\langle v, r_j \right\rangle|^2 } ~ \mbox{satisfies} ~ \left(1 - \varepsilon\right)\frac{m}{n} \leq h(v) \leq \left(1 + \varepsilon\right)\frac{m}{n} \quad \mbox{w.h.p.}$$
This is, in fact, exactly the question that underlies the study of random projections, and has been dealt with extensively (see e.g., \cite{vempala}). The main conclusion is that
$m \sim \log{n}$ is in many cases sufficient, however, the implicit constant may be large.\\ 

We now explain why in practice a much smaller number of random vectors
suffices. The idea is rather simple and most easily explained by
considering the example of one random vector in the case of a large
spectral gap $\lambda_1 \gg \lambda_2$. Clearly, we have
\begin{equation*}
 A^{\alpha} r = \sum_{i=1}^{n}{ \lambda_i^{\alpha} \phi_i \left\langle \phi_i, r \right\rangle} =  \lambda_1^{\alpha} \phi_1 \left\langle \phi_1, r\right\rangle +  \sum_{i=2}^{n}{ \lambda_i^{ \alpha} \phi_{i} \left\langle \phi_i, r \right\rangle}.
\end{equation*}
The outcome now depends on the random vector $r$, however, for
suitable large values of $\alpha$ it is clear that in order for the
random landscape to profoundly differ from the profile of the leading
eigenvector, it is required that
$|\left\langle \phi_1, r \right\rangle|$ is very small: even if it
were only moderately small, it would get drastically amplified by the
exponential growth and still dominate the expression. The following
widely-used Lemma shows that this is not overly likely.

\begin{lemma} Let $v \in \mathbb{R}^n$ satisfy $\|v\| = 1$ and let $r \in \mathbb{S}^{d-1}$
  be randomly chosen w.r.t. the uniform surface measure on
  $\mathbb{S}^{n-1}$. Then, for $0 \leq \delta \leq 1$,
  $$ \mathbb{P}\left( |\left\langle v, r \right\rangle| \leq \frac{\delta}{\sqrt{n}} \right) \lesssim \delta,$$
  where the implicit constant depends only on the dimension.
\end{lemma}

This simple Lemma quantifies the concentration of measure phenomenon and is
standard (see e.g. \cite{ken}). It explains why in the case of highly localized
eigenstates it is completely sufficient to work with only one randomly
chosen vector. The inner products
$\left\langle \phi_i, r \right\rangle$ are not likely to be extremely
small and get amplified by an exponential growth while the strong
exponential localization preserves the structure. This is easy to make
precise in a variety of ways: the simplest case is a spectral gap $\lambda_1 \gg \lambda_2$.
The trivial estimate
$$ \left\| \sum_{i=2}^{n}{ \lambda_i^{\alpha} \phi_{i} \left\langle \phi_i, r \right\rangle} \right\|_{\ell^2} \leq   |\lambda_2|^{\alpha}\left(\sum_{i=2}^{n}{ \left\langle \phi_i, r \right\rangle^2}\right)^{1/2} \leq  |\lambda_2|^{\alpha} $$
implies, together with the Lemma above,
$$ \mathbb{P}\left( \left\|  \lambda_1^{\alpha} \phi_{1} \left\langle \phi_1, r\right\rangle \right\|_{\ell^2} \leq \left\| \sum_{i=2}^{n}{ \lambda_i^{\alpha} \phi_{i} \left\langle \phi_i, r \right\rangle} \right\|_{\ell^2} \right) \lesssim \sqrt{n} \left( \frac{|\lambda_2|}{|\lambda_1|}\right)^{\alpha}.$$
We emphasize that this simple argument did not even use localization of the eigenvectors; the trivial estimate is clearly quite weak if the spectrum is spread out, in that case much stronger results should hold.
\medskip 

\textit{Example.} Let us revisit the example in
Section~\ref{sec:mainidea} using the random sampling.

\begin{figure}[ht]
\centering
\begin{subfigure}[b]{0.49\textwidth}
\includegraphics[width = 0.95\textwidth]{bandmatland.pdf} 
\caption{$f_{50}$ without random sampling}
\end{subfigure}
\begin{subfigure}[b]{0.49\textwidth}
\includegraphics[width = 0.95\textwidth]{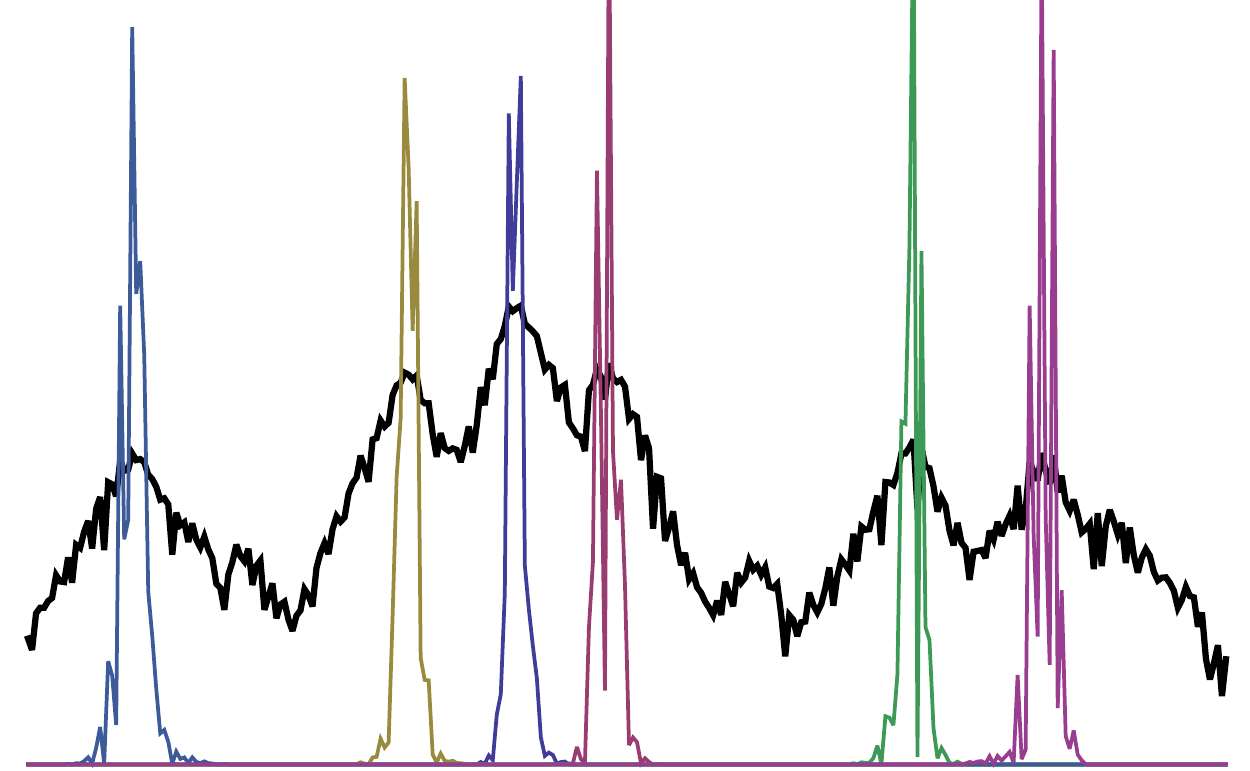} 
\caption{$f_{R, 50}$ with $1$ random vector}
\end{subfigure} \\
\begin{subfigure}[b]{0.49\textwidth}
\includegraphics[width = 0.95\textwidth]{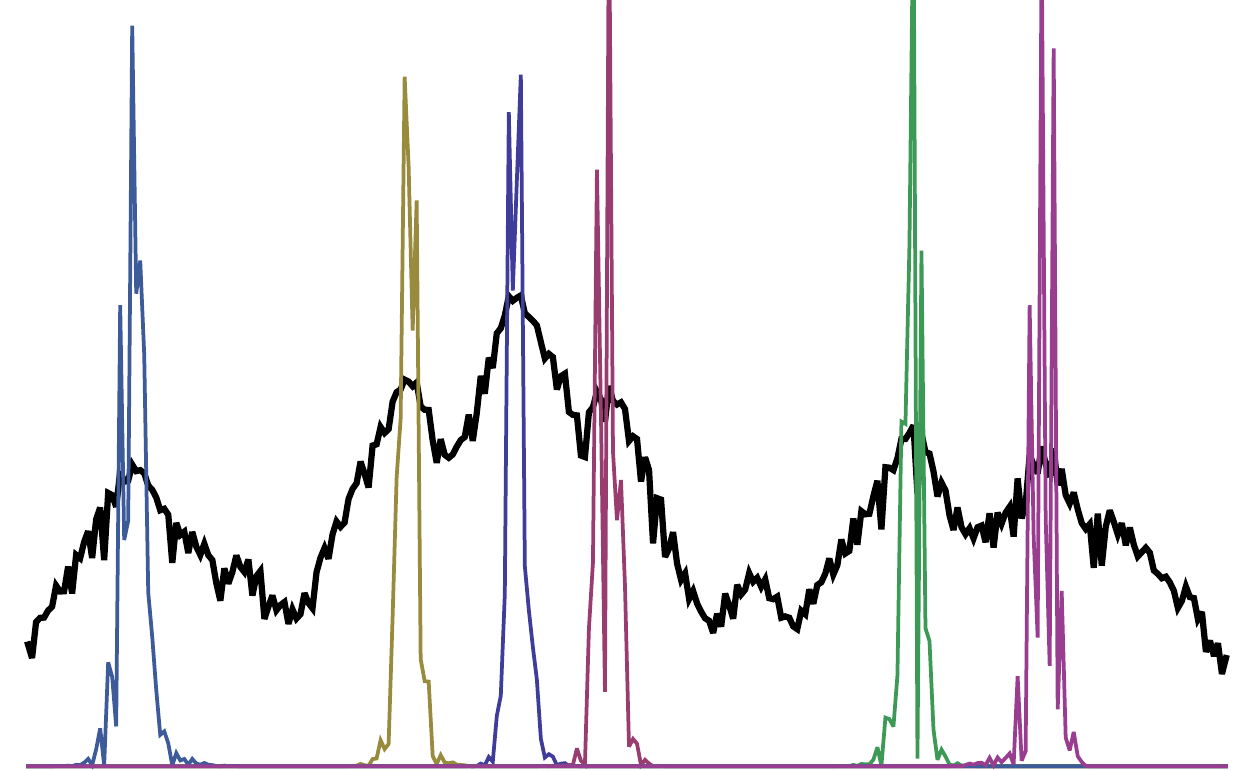} 
\caption{$f_{R, 50}$ with $3$ random vectors}
\end{subfigure} 
\begin{subfigure}[b]{0.49\textwidth}
\includegraphics[width = 0.95\textwidth]{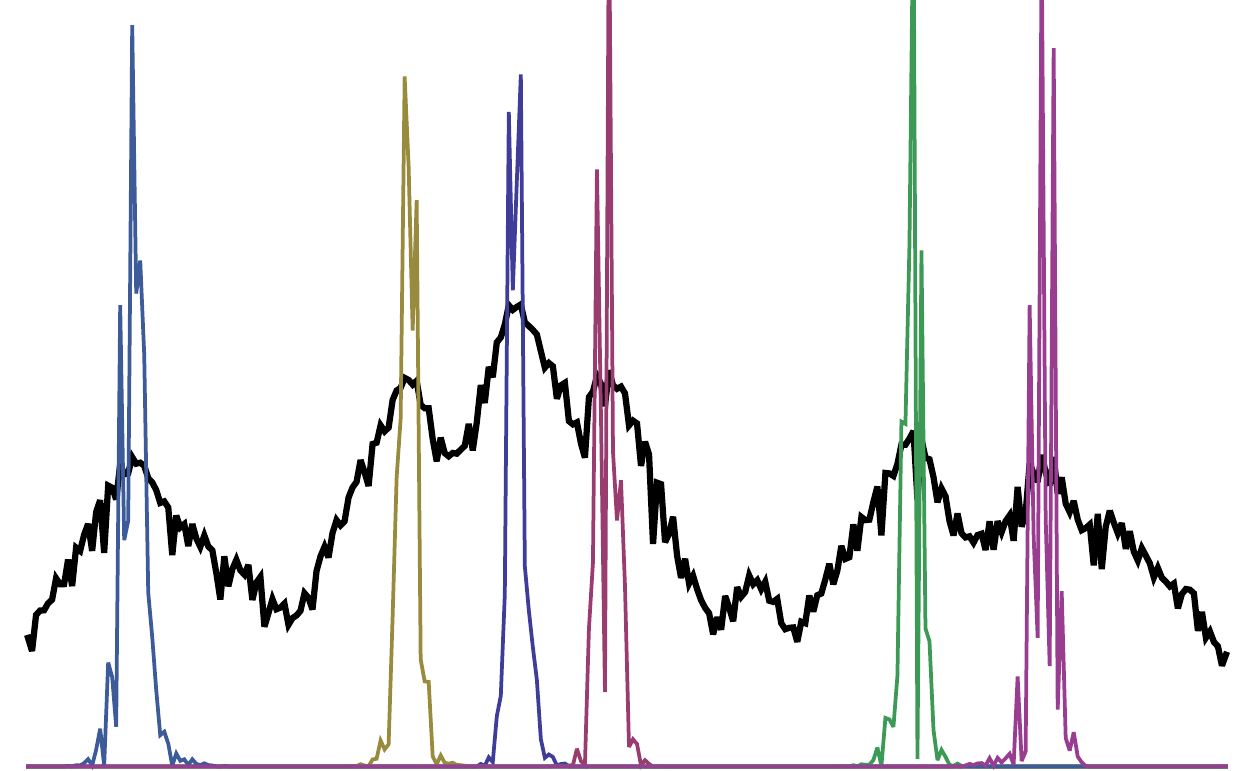} 
\caption{$f_{R, 50}$ with $5$ random vectors}
\end{subfigure}
\caption{The landscape function $f_{50}$ (A), compared
  with its randomized sampled version with $1$ (B), $3$
  (C), and $5$ (D) random vectors, are plotted
  in black.  The absolute value of the first $6$ eigenvectors
  associated to the $6$ largest eigenvalues is also plotted. 
   (A) is a replication of Figure~\ref{fig:eins}.}
\label{fig:eins_random}
\end{figure}

 The same random
band matrix $A \in \mathbb{R}^{n \times n}$, $n=300$ for
Figure~\ref{fig:eins} is used. The randomized version with number of
random vectors $1$, $3$, $5$ respectively are plotted in the panels of
Figure~\ref{fig:eins_random}. We observe that the randomized landscape
function even with only $1$ random vector still captures the important
feature, in particular the local maxima, of the landscape function.

\section{Numerical Examples} \label{sec:num}

\subsection{Schr\"odinger Operator with Potential.} 
The case of finding a way of numerically detecting localized
eigenstates of operators of the form $-\Delta + V$ for $V \geq 0$ and
Dirichlet conditions on the boundary has recently received a lot
attention in the mathematics and
physics literature (see e.g., \cite{ lots, fil,fil2, fil3, stein}). We explain how our method can be applied to this
case (without any restrictions on $V$).
We consider the operator $-\Delta + V$ defined on $[0, 1]^2$ with
periodic boundary condition. $V$ is a smooth periodic potential
generated randomly in the unit square. The pseudo-spectral
discretization with mesh size $h = 1/48$ is used. 
\begin{figure}[ht]
\centering
\begin{subfigure}[b]{0.49\textwidth}
\includegraphics[width = 0.85\textwidth]{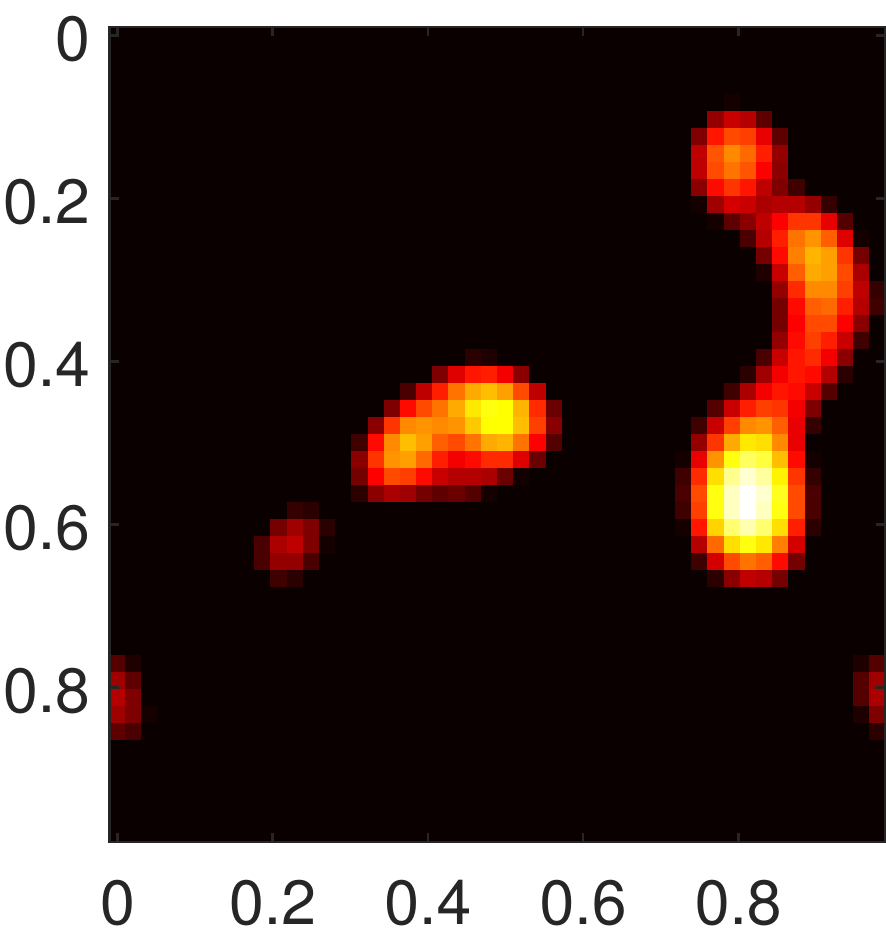} 
\caption{$f_{R, 75}$ with $10$ random vectors}
\end{subfigure}
\begin{subfigure}[b]{0.49\textwidth}
\includegraphics[width = 0.85\textwidth]{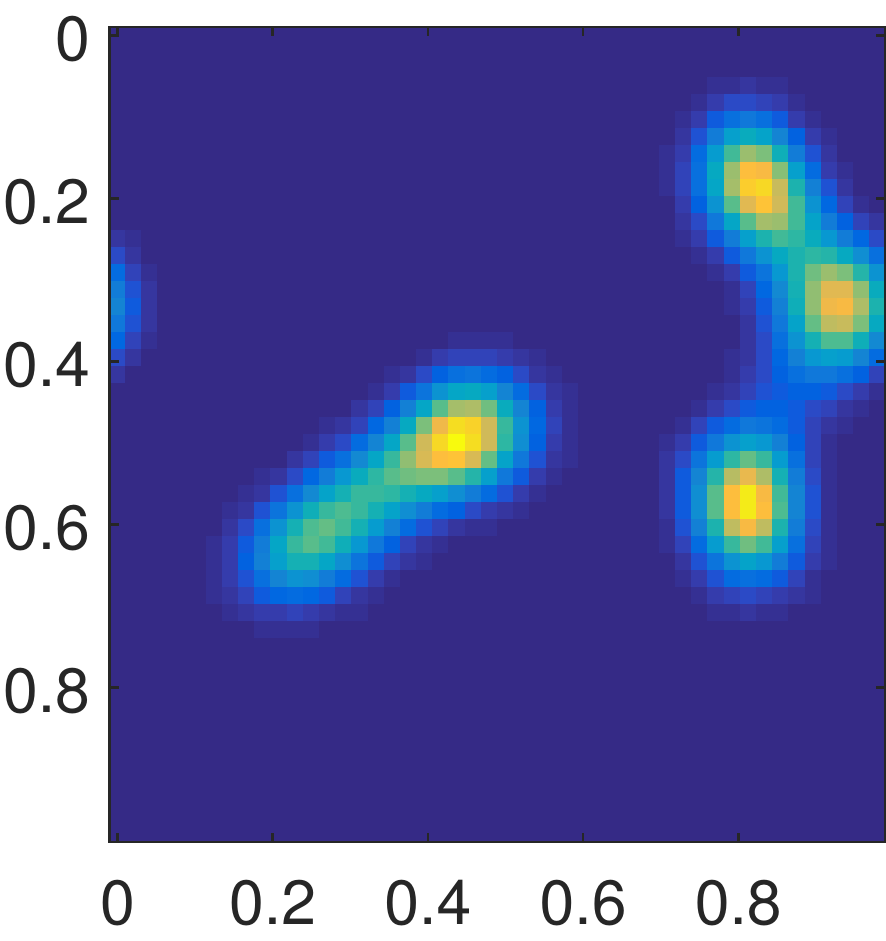} 
\caption{Sum of modulus square:
  $\sum_{i=1}^5 |\phi_i|^2$ }
\end{subfigure} \\
\begin{subfigure}[b]{0.49\textwidth}
\includegraphics[width = 0.85\textwidth]{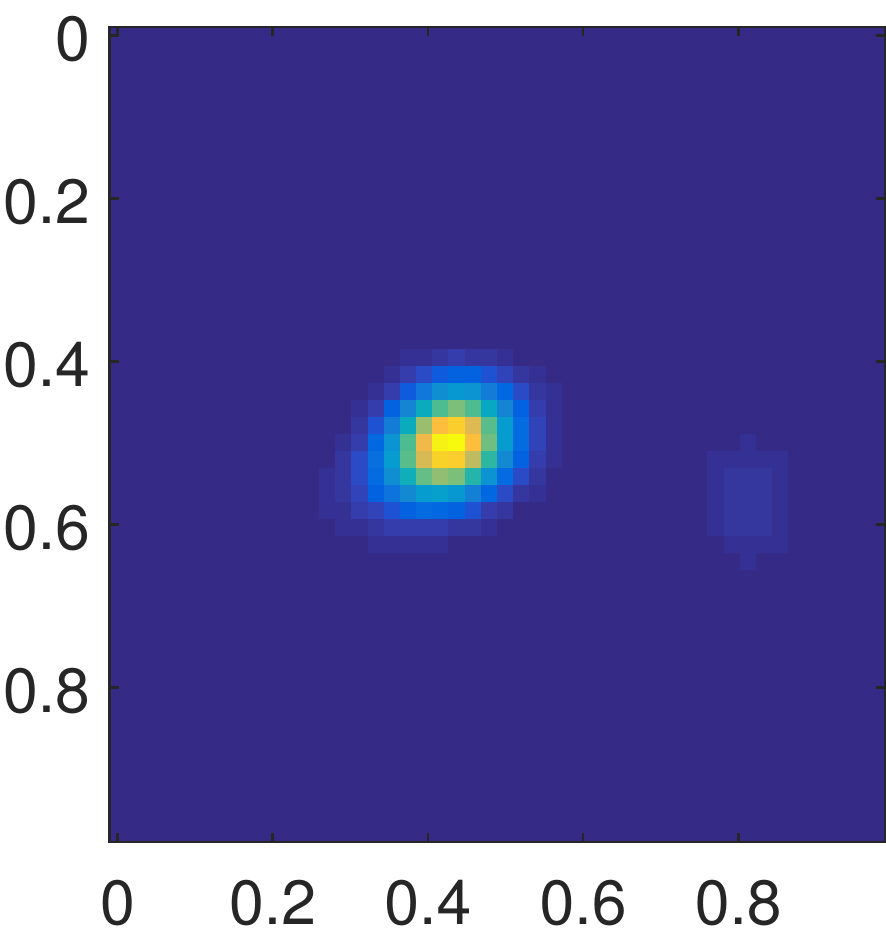} 
\caption{First low-lying eigenfunction $|\phi_1|^2$}
\end{subfigure} 
\begin{subfigure}[b]{0.49\textwidth}
\includegraphics[width = 0.85\textwidth]{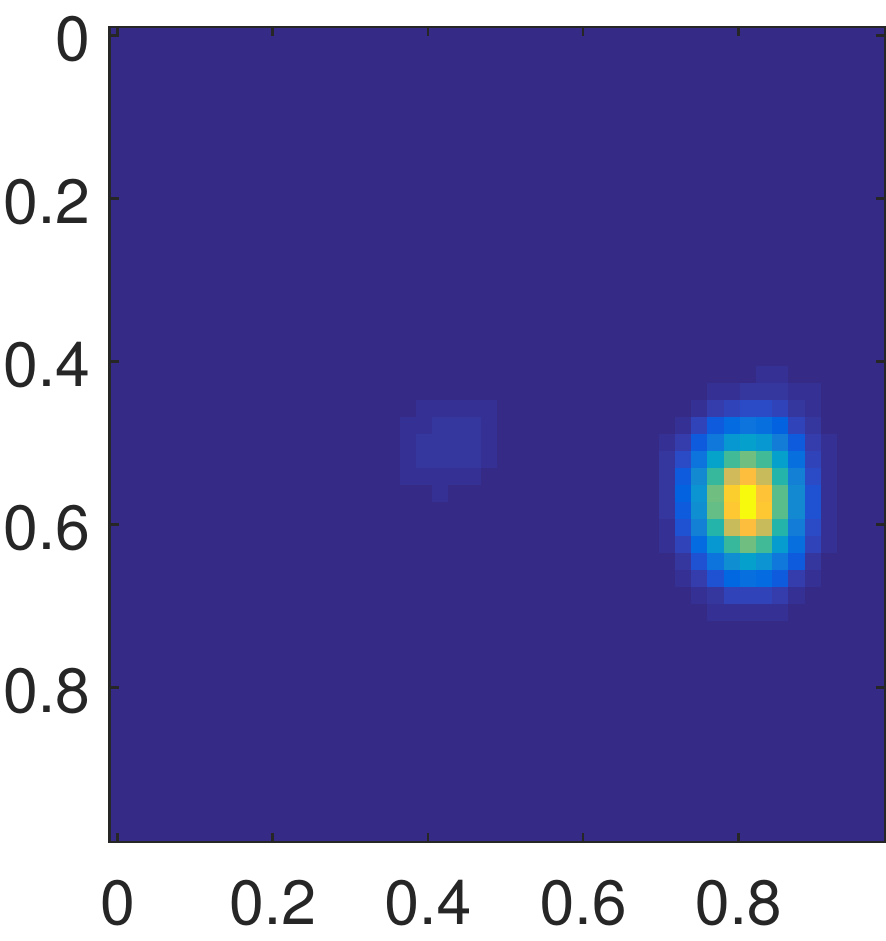} 
\caption{Second low-lying eigenfunction $|\phi_2|^2$}
\end{subfigure}
\caption{The landscape function $f_{R, 75}$ (A), compared
  with the sum of modulus square of the low-lying $5$ eigenfunctions (B), for the operator $-\Delta + V$. The modulus square of the first two low-lying eigenfunctions are also plotted in (C) and (D) respectively. }
\label{fig:schr}
\end{figure}
While the
resulting matrix is dense due to the Fourier differentiation of the
pseudo-spectral method, the method proposed still applies. In
Figure~\ref{fig:schr}, we show the landscape function using
$\alpha = 75$ and $10$ random vectors on the left panel and the sum of
the square modulus of the low-lying $5$ eigenvectors on the right
panel. Good agreement is observed.  For visualization, we only plotted
the part of landscape function exceeding its maximum value minus $1$.

\begin{figure}[ht]
\centering
\begin{subfigure}[b]{0.49\textwidth}
\includegraphics[width = 0.95\textwidth]{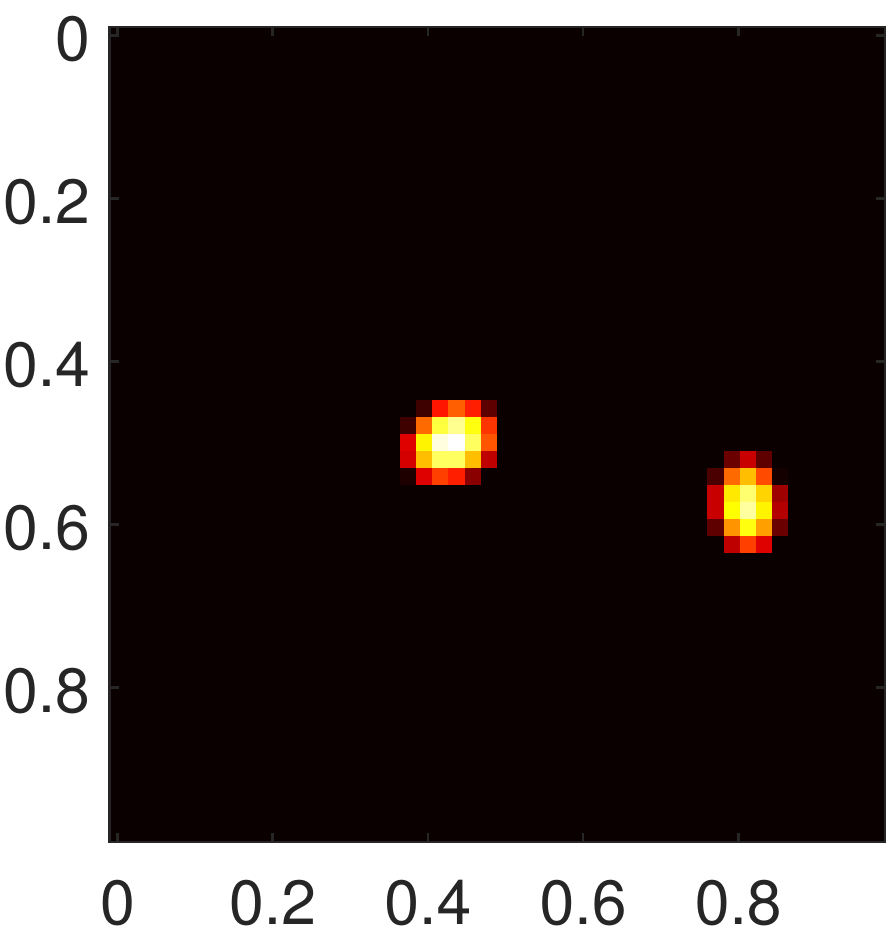} 
\caption{$f_{R, 75}$ with $10$ random vectors}
\end{subfigure}
\begin{subfigure}[b]{0.49\textwidth}
\includegraphics[width = 0.95\textwidth]{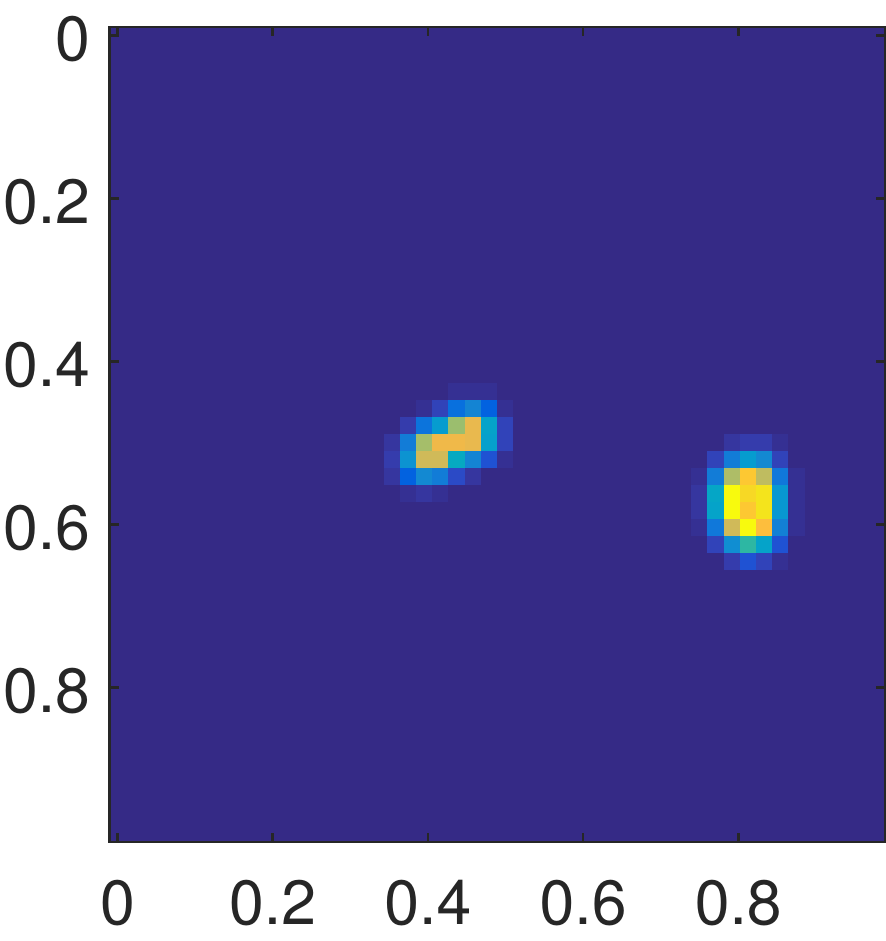} 
\caption{Sum of modulus square:
  $\sum_{i=1}^5 |\phi_i|^2$ }
\end{subfigure} \\
\begin{subfigure}[b]{0.49\textwidth}
\includegraphics[width = 0.95\textwidth]{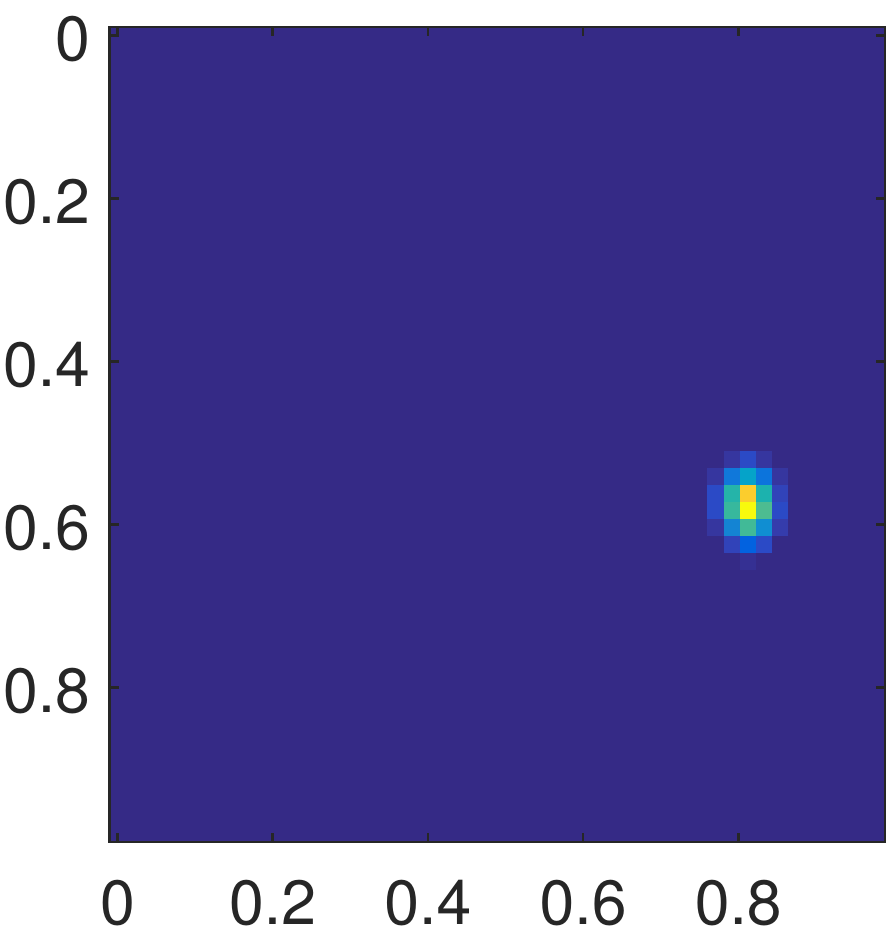} 
\caption{First low-lying eigenfunction $|\phi_1|^2$}
\end{subfigure} 
\begin{subfigure}[b]{0.49\textwidth}
\includegraphics[width = 0.95\textwidth]{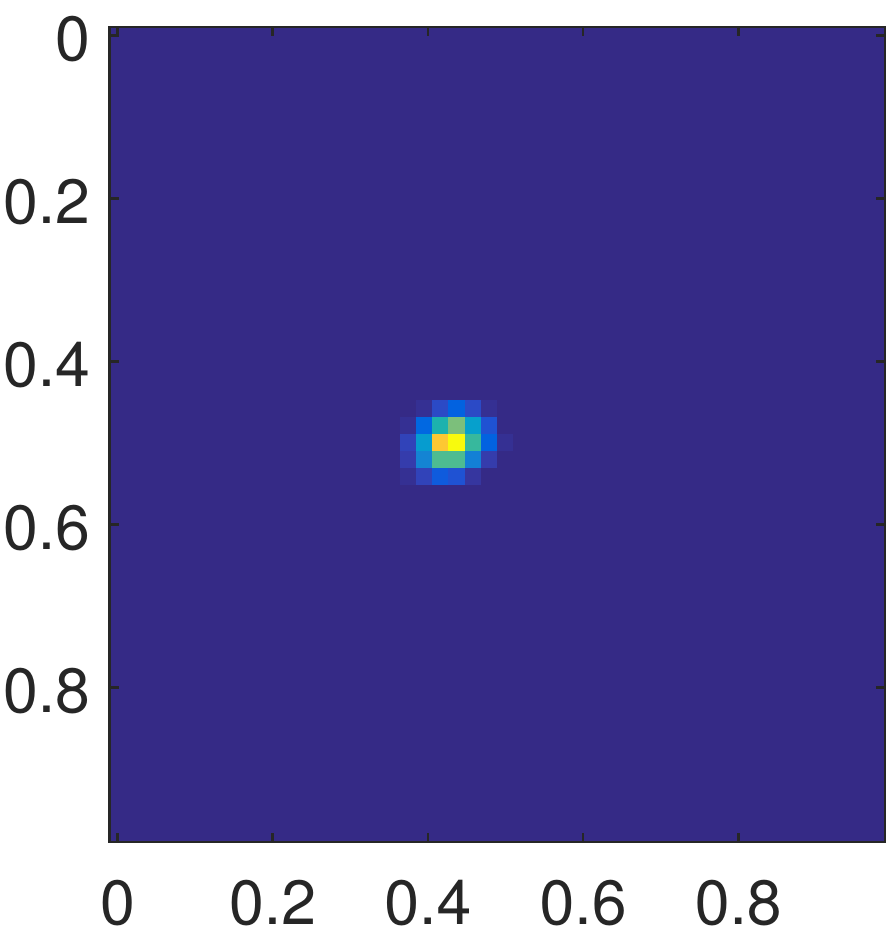} 
\caption{Second low-lying eigenfunction $|\phi_2|^2$}
\end{subfigure}
\caption{The landscape function $f_{R, 75}$ (A), compared
  with the sum of modulus square of the low-lying $5$ eigenfunctions
  (B), for the operator $(-\Delta)^{3/4} + V$. The same
  potential as in the example of Figure~\ref{fig:schr} is used. The
  modulus square of the first two low-lying eigenfunctions are also
  plotted in (C) and (D) respectively. }
\label{fig:frac34}
\end{figure}

\subsection{Fractional Laplacian}

The method extends without much difficulty to the fractional Laplacian
$(-\Delta)^{s} + V$. We still consider the computational domain
$[0, 1]^2$ with periodic boundary condition and $V$ a smooth randomly generated potential. 
Thanks to the periodic boundary condition, the
fractional Laplacian $(-\Delta)^s$ can be defined through spectral
decomposition. Pseudo-spectral discretization is again used. Note that
the fractional Laplacian is non-local, regardless of the
discretization. In Figure~\ref{fig:frac34}, we show the result for
$(-\Delta)^{3/4} + V$ where the same potential as in the example of
Figure~\ref{fig:schr} is used. We again observe excellent agreement
between the landscape function and the localization of the
eigenfunctions.

\section{Proofs}
\begin{proof}[Proof of the Theorem] We first show that $f_{\alpha}$ exceeds a certain value on all the intervals $J_i$ ($1 \leq i \leq k$). We start by showing there is a $k \in J_i$ with
$$  \left| \left\langle \phi_i, e_k \right\rangle \right| \geq \frac{1}{2} \frac{1}{\sqrt{\mathcal{J}}}.$$
If this were false, then
$$ \frac{1}{4} \leq   \sum_{k \in J_i}{ \left| \left\langle \phi_i, e_k \right\rangle \right|^2} <  \sum_{k \in J_i}{  \frac{1}{4J}} = \frac{|J_i|}{4 \mathcal{J}} \leq \frac{1}{4},$$
which is a contradiction. Here the first inequality comes from the
assumption that $J_i$ contains half $\ell^2$-mass of $\phi_i$.  A
simple application of the spectral theorem implies that
$$ \| A^{\alpha} e_k\|_{\ell^2} \geq \frac{1}{2\sqrt{\mathcal{J}}} |\lambda_i|^{\alpha} \geq  \frac{1}{2\sqrt{\mathcal{J}}} |\lambda_k|^{\alpha}.$$
 It remains to show that this value cannot be attained unless one is close to one of the localized eigenvectors.
If $m \in \left\{1,2, \dots, n\right\}$ has at least distance $d>0$ from $\bigcup_{i=1}^{k}{J_i},$ then 
\begin{align*}
 \| A^{\alpha} e_m\|^2_{\ell^2} &= \sum_{l=1}^{n}{ \left|\left\langle e_m, \phi_l \right\rangle \right|^{2} |\lambda_l|^{2\alpha}} \leq e^{-2 \beta d} k |\lambda_1|^{2\alpha} + n |\lambda_{k+1}|^{2\alpha}
\end{align*}
and thus
$$   \| A^{\alpha} e_m\|_{\ell^2} \leq  e^{- \beta d}\sqrt{k} |\lambda_1|^{\alpha} + \sqrt{n} |\lambda_{k+1}|^{\alpha}.$$
We want to relate this to the inequality
$$ e^{-\beta d} \sqrt{k} |\lambda_1|^{\alpha} + \sqrt{n} \left(|\lambda_{k+1}| \right)^{\alpha} \leq \frac{1}{4} \left(   \frac{1}{2\sqrt{\mathcal{J}}} |\lambda_k|^{\alpha} \right),$$
which will provide an upper bound on $d$ (meaning that if the distance is to the sets is larger than this bound, then we obtain the desired consequence).
If, as we assume,
$$\left( \frac{\lvert \lambda_k \rvert}{\lvert \lambda_{k+1} \rvert} \right)^{\alpha} \geq 16 \sqrt{\mathcal{J}n} \qquad \mbox{then} \qquad  \sqrt{n} \left(|\lambda_{k+1}| \right)^{\alpha} \leq  \frac{1}{8} \left(   \frac{1}{2\sqrt{\mathcal{J}}} |\lambda_k|^{\alpha} \right).$$
The inequality
$$  e^{-\beta d} \sqrt{k} |\lambda_1|^{\alpha}  \leq \frac{1}{8} \left(   \frac{1}{2\sqrt{\mathcal{J}}} |\lambda_k|^{\alpha} \right)$$
is equivalent to
$$ d \leq  \frac{1}{ \beta} \left(\log{\left(16\sqrt{k}\sqrt{\mathcal{J}}\right)} + \alpha \log{\left( \left| \frac{\lambda_1}{\lambda_k} \right| \right)}\right)$$
which is the desired statement.
\end{proof}


\appendix

\section{Localized Wannier Bases}

In this appendix, we connect the landscape functions to localized
Wannier bases.  Localized Wannier bases \cite{wannier} are maximally
localized bases of low-lying eigenfunctions. Give an operator over
$\Omega$ with mutually orthogonal low-lying eigenfunctions
$\phi_1, \dots, \phi_n \in L^2(\Omega)$, it is often helpful to work
with a basis of $\mbox{span}\left\{\phi_1, \dots, \phi_n \right\}$
where each basis element is as localized as possible (for example to
obtain sparser matrices).
A classical approach, see e.g., the review article \cite{Marzari:12}, is to simply project the Dirac measure $\delta_x$ onto the finite-dimensional subspace $\mbox{span}\left\{\phi_1, \dots, \phi_n \right\}$ to obtain the best possible representation of this localized point in the basis. We illustrate this with a simple
example: let us consider $\left\{e^{-i n x}, e^{-i (n-1) x}, \dots, 1, \dots, e^{i n x}\right\}$ on the torus $\mathbb{T}$. The
projection of $\delta_y$ onto the span is given by
$$ \pi_n \delta_y = \sum_{k}{\left\langle \delta_y, \phi_k\right\rangle \phi_k} = \sum_{k=-n}^{n}{ e^{i k y} e^{-i k x}} =  \sum_{k=-n}^{n}{ e^{i k (x+y)}} = \frac{ \sin{\left(\left(n + \frac{1}{2}\right) (x-y)\right)}}{\sin{(x-y)/2}},$$
which is merely the classical Dirichlet kernel. The function is indeed highly localized around $y$, roughly constant for $|x-y| \lesssim n^{-1}$ and then exhibits some oscillatory behavior around 0 further away from  $y$. Motivated by our main idea, a slightly different idea suggests itself: instead of merely projecting, it makes sense to run the dynamical system with $\delta_y$ as initial datum on the subspace on which we project. More precisely, this yields a projection that depends on $t$
$$ \pi_{n,t} \delta y :=   \sum_{k}{e^{-\lambda_k t}\left\langle \delta_y, \phi_k\right\rangle \phi_k}.$$

\begin{center}
\begin{figure}[h!]
\includegraphics[width=0.7\textwidth]{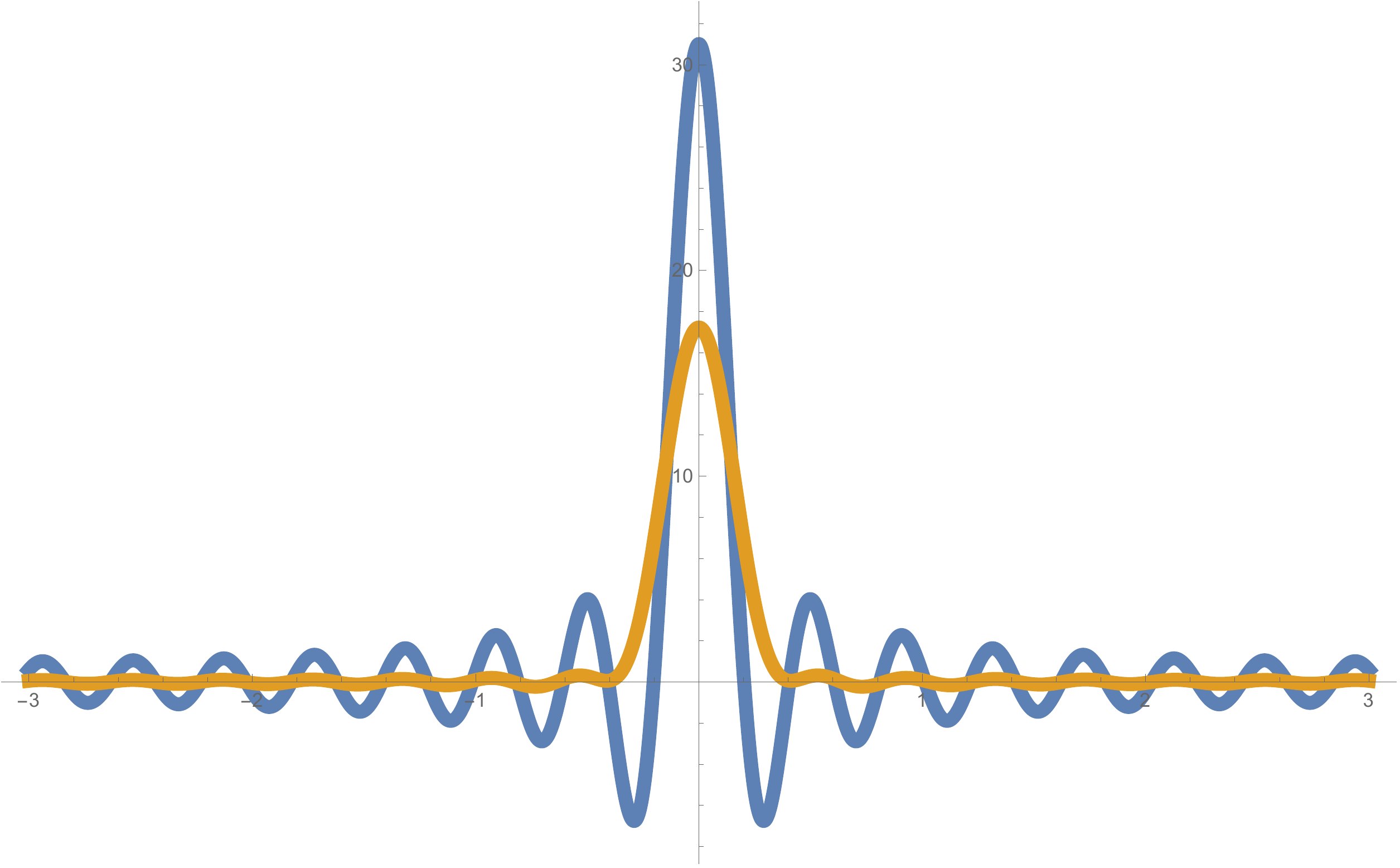}
\caption{The projection on $\mathbb{T}$ for $n=15$ and the parametrized projection with $t=0.01$. We observe a tradeoff between localization and decay.}
\end{figure}
\end{center}

Basic intuition tells us that, depending on the speed of propagation within the dynamical system, at least for small values of $t$ the projection will still look pretty much like the direct propagation. Our main observation is that for sufficiently small values of $t$, the
projection $\pi_{n,t}$ will be pretty much as localized as $\pi_n$ but will have much better decay properties. We illustrate this again on the torus $\mathbb{T}$, where
$$ pi_{n,t} \delta y :=   \sum_{k}{e^{-\lambda_k t}\left\langle \delta_y, \phi_k\right\rangle \phi_k} = \sum_{k}{e^{-k^2 t} e^{-iky} e^{i kx}}.$$
We observe that the infinite limit is given by the Jacobi $\theta-$function
$$ \theta_t(x-y) = \sum_{k = \infty}^{\infty}{ e^{-k^2 t} e^{i k (x-y)}}.$$
The Jacobi $\theta-$function has exponential decay in the region $|x-y| \gtrsim t^{1/2}$ and, what is especially useful, the exponential decay in the coefficients implies that a cutoff at frequency $n$ is only going to introduce a small error. In particular,
in this setting, the threshold-value for $t$ (where localization properties are not significantly worse) is $t \lesssim n^{-2}$. This example generalizes immediately to higher-dimensions, where the situation is more or less identical: the direct projection
$\pi_n$ creates oscillations with slow decay throughout the entire space, the diffused $\pi_{n,t}$ leads to a more localized representation. 
Returning to our original problem for matrices $A \in \mathbb{R}^{n \times n}$, the relationship to $f_{\alpha}$ is easily explained: in the case of localized eigenvectors, we  see that if $\phi$ is the eigenvector associated to the largest eigenvalue on coordinate $k \in \left\{1, 2, \dots, n\right\}$, then $A^{\alpha} \left( \delta_{ik}\right)_{i=1}^{n}$ is a very good local approximation of $\phi$.
The proof of this simple statement proceeds along the very same lines as the proof of our main theorem.

\end{document}